\documentstyle[12pt, amscd]{amsart}

\newtheorem{theorem}{Theorem}[section]
\newtheorem{lemma}[theorem]{Lemma}

\theoremstyle{definition}
\newtheorem{definition}[theorem]{Definition}

\theoremstyle{remark}

\numberwithin{equation}{section}

\newcommand{\HE}{\mathrm {HE}}
\newcommand{\Mat}{\mathrm {Mat}}
\newcommand{\HP}{\mathrm {HP}}

\newcommand{\Lie}{\mathrm {Lie}}
\newcommand{\GL}{\mathrm {GL}}
\newcommand{\SL}{\mathrm {SL}}

\newcommand{\Hom}{\mathrm {Hom}}
\newcommand{\ad}{\mathrm {ad}}
\newcommand{\Tot}{\mathrm {Tot}}

\begin{document}

\title{Non-Commutative Chern Characters of Compact Quantum Groups}

\author{Do Ngoc Diep} 
\address{Institute of Mathematics, National
Centre for Natural Science and Technology\\ P.O. Box 631, Bo Ho, 10.000,
Hanoi, Vietnam} 
\curraddr{International Centre for
Theoretical Physics, ICTP P. O. Box 586, 34100, Trieste, Italy}
\email{dndiep@@ioit.ncst.ac.vn} 
\thanks{The first author was supported in part by the International Centre
for Theoretical Physics, the IMU Commission on Development and Exchange,
and Vietnam National Foundation for Research in Fundamental Sciences}

\author{Aderemi O. Kuku}
\address{International centre for Theoretical Physics, ICTP P. O. Box 586, 34100,
Trieste, Italy}
\email{kuku@@ictp.trieste.it}

\author{Nguyen Quoc Tho}
\address{Department of Mathematics, Vinh University, Vinh City, Vietnam}
\curraddr{c/o Institute of Mathematics, National Centre for Natural Science and
Technology\\
P.O. Box 631, Bo Ho, 10.000, Hanoi, Vietnam}
\email{dndiep@@ioit.ncst.ac.vn}

\date{July 17, 1998} 


\begin{abstract}
In this paper we compute the K-theory (algebraic and topological) and
entire periodic cyclic homology for compact quantum groups, define Chern
characters between them and show that the Chern characters in both
topological and algebraic cases are isomorphisms modulo torsion.
\end{abstract}
\maketitle

\section*{Introduction}

Non commutative differential geometry appeared as some quantization of
classical differential geometry. In particular, different cyclic
(co)homologies are non-commutative versions of the classical de Rham
(co)homology of (differential forms or) currents. Quantum groups were
introduced, at first as Hopf bialgebra envelop of Lie groups and later as
special cases of quasi-triangular Drinfeld algebras or monoidal categories
. The theory is a subject of intensive current research. 

Many results have
been obtained for quantum groups as analogues of the corresponding
classical results for Lie groups, especially the theory of highest weight
representations. For compact quantum groups, the representation theory is
characterized by the highest weight representations, just as in the case
of representation theory of compact Lie groups. For the compact Lie
groups, the $K$-groups, de Rham cohomology groups and the Chern
characters between them have been well studied (see \cite{H},
\cite{W}-\cite{W'}). More recently the Chern characters from K-theory
to entire cyclic homology of such groups have also been studied (see
\cite{DKT}). However  such studies have not yet been done for
quantum groups. The aim of this paper is to compute the non-commutative
Chern characters for compact quantum groups in a C*-algebra setting. 

Algebraic and topological $K$-theories of Banach algebras have been well
studied and known to have such nice properties as homotopy invariance,
Morita invariance, excision, etc... At the same time, the entire cyclic
cohomology $\HE^*(A)$, defined by A. Connes with the non-commutative Chern
characters $ch_{C^*}$, given by the pairing $ K_*(A)  \times \HE^*(A) \to
{\mathbf C}$ has been rather difficult to compute for various concrete
algebras and indeed the authors are not aware of any computations
involving algebras other
then $A ={\mathbf C}$. 

In this paper we define an entire cyclic homology $\HE_*(A)$ by
considering a direct system of ideals of $A$ with $ad$-invariant
trace and using the pairing of the entire cochains of A. Connes at
the level of ideals with elements of $K_*(A)$

We now the structure of the paper. In Section 1, we at first
recall definitions and main properties of compact quantum groups and their
representations. More precisely, if $G$ is a complex algebraic group with
compact real form, ${\mathfrak g}$ its Lie algebra, $\varepsilon$ a
positive real number, we define a C*-algebraic compact quantum group
$C^*_\varepsilon(G)$ as the C*-completion of the *-algebra ${\mathcal
F}_\varepsilon(G)$ with respect to the C*-norm, where ${\mathcal
F}_\varepsilon(G)$ is the quantized Hopf algebra of functions which is
the Hopf subalgebra of the Hopf algebra dual to the quantized universal
enveloping algebra $U({\mathfrak g})$, generated by matrix elements of
the  $U({\mathfrak g})$ modules of type ${\mathbf 1}$(see \cite{CP}). We
proved that 
$$C^*_\varepsilon(G) \cong {\mathbf C}({\mathbf T}) \oplus \bigoplus_{e
\ne w\in W}\int_{\mathbf T}^{\oplus} {\mathcal K}({\mathbf H}_{w,t}) dt$$
where ${\mathcal K}({\mathbf H}_{w,t})$ is the elementary
algebra of compact operators in separable
infinite-dimensional Hilbert space ${\mathbf H}_{w,t}$ and $W$ is the
Weyl group of $G$ with respect to a maximal torus ${\mathbf T}$. 

In Section 2, we provide the definition and main properties of entire
periodic cyclic homology of non-commutative involutive Banach algebras. We
prove in particular that when the C*-algebra $A$ can be presented as
direct limit $\varinjlim I_\alpha$, where $I_\alpha$ are ideals of $A$
admitting an $\ad$-invariant trace then $\HE_*(A) = \varinjlim
HE_*(I_\alpha)$. 

Section 3 is devoted to Chern characters. We first compute the
K-groups of $C^*_\varepsilon(G)$, the $\HE_*(C^*_\varepsilon(G))$.
Therefore we  define the Chern characters from $K_*(C^*_\varepsilon(G))$
to $\HE_*(C^*_\varepsilon(G))$ for any involutive Banach algebra $A$
endowed with a family of ideals with $\ad$-invariant trace. When
$A=C^*_\varepsilon(G)$ we prove that
$$\HE_*(C^*_\varepsilon(G)) \cong \HE_*({\mathbf C}({\mathcal N}_{\mathbf
T})) \cong H^*_{DR}({\mathcal N}_{\mathbf T})$$ and
$$K_*(C^*_\varepsilon(G)) \cong K_*({\mathbf C}({\mathcal N}_{\mathbf
T})) \cong K^*({\mathcal N}_{\mathbf T}).$$ 

Finally we show that there is a commutative diagram 
$$\CD
K_* (C^*_\varepsilon(G)) @>ch_{C^*}>> \HE_*(C^*_\varepsilon(G))\\
@V\cong VV    @VV\cong V\\
K_* ({\mathbf C}({\mathcal N}_{\mathbf T})) @>ch_{CQ}>>
\HE_*({\mathbf C}({\mathcal N}_{\mathbf T}))\\
@V\cong VV    @VV\cong V\\
K^* ({\mathcal N}_{\mathbf T}) @>ch> > H^*_{DR}({\mathcal
N}_{\mathbf T}),
\endCD$$ 
we deduce that $ch_{C^*}$ is an isomorphism modulo torsion.

In Section 4 we at first identify our definition of
$\HP_*(C^*_\varepsilon(G))$ with Cuntz-Quillen definition, and then
compute
$K_*^{alg}(C^*_\varepsilon(G))$, $\HP(C^*_\varepsilon(G))$ and prove that
the Chern character $ch_{alg}$ from $K_{alg}(C^*_\varepsilon(G))$ to
$\HP_*(C^*_\varepsilon(G))$ is an isomorphism modulo torsion.

{\sl Notes on Notation}. For any compact space $X$, we write $K^*(X)$ for
the ${\mathbf Z}/2$ graded topological K-theory of $X$. We use Swan's
theorem to identify $K^*(X)$ with ${\mathbf Z}/2$ graded $K_*({\mathbf
C}(X))$. For any involutive Banach algebra $A$, $K_*(A)$, $\HE_*(A)$,
$\HP_*(A)$ are ${\mathbf Z}/2$ graded algebraic or topological K-group of
$A$, entire cyclic homology, and periodic cyclic homology of $A$,
respectively. If ${\mathbf T}$ is a maximal torus of a compact group $G$,
with the corresponding Weyl group $W$, write ${\mathbf C}({\mathbf T})$
for the algebra of complex valued functions on ${\mathbf T}$. We use the
standard notations from the root theory such as $P$, $P^+$ for the
positive highest weights, etc....  We denote by ${\mathcal N}_{\mathbf T}$
the normalizer of ${\mathbf T}$ in $G$, by ${\mathbf N}$ the set of
natural numbers, ${\mathbf R}$ the field of real numbers and ${\mathbf
C}$
the field of complex numbers, $\ell^2_A({\mathbf N})$ the standard
$\ell^2$ space of square integrable sequences of elements from $A$, and
finally by $C^*_\varepsilon(G)$ we denote the compact quantum group. 

\section{Compact Quantum Groups and Their Representations}

In this section we present a summary of some notions and results
concerning compact quantum groups and their representations for later use.
The main reference is the book \cite{CP}

Let $G$ be a complex algebraic group with Lie algebra ${\mathfrak g} =
\Lie\; G$. There is a one-to-one correspondence between:
\begin{enumerate}
\item real forms of $G$,
\item Hopf *-structures on the universal enveloping algebra $U({\mathfrak
g})$ and
\item Hopf *-structures on the dual ${\mathcal F}(G) := U({\mathfrak g})'$.
\end{enumerate}

The compact real form $K$ can be characterized as the set of fixed
points of conjugate-linear involution of the algebra ${\mathcal F}(G)$ of
representative functions of $G$. 

In the Hopf algebra dual $U_\varepsilon({\mathfrak g})'$, one considers
the Hopf sub-algebra ${\mathcal F}_\varepsilon(G)$, generated by the matrix
elements of finite-dimensional $U_\varepsilon({\mathfrak g})$-module of
type ${\mathbf 1}$, see Definition 13.1.1 from \cite{CP}.

If $\varepsilon$ is real and $\varepsilon \ne -1$, there exists a unique Haar
functional on ${\mathcal F}_\varepsilon(G)$ which can be used to prove results
similar to those for classical Lie groups, e.g. the Schur orthogonality relations.
From this one can effectively study the unitary representations of compact quantum
groups, see (\cite{CP}, Ch. 13). 

In the classical theory, there is a (one-dimensional) unitary representation of the
algebra ${\mathcal F}(G)$ associated to every point of the real form $G_{\mathbf R}$
the corresponding involution.  Since in the quantum case there are no ``points", the
study of unitary representations replaces the study of the space itself. For compact
quantum groups the unitary representations of ${\mathcal F}_\varepsilon(G)$ are
parameterized by pairs $(w,t)$, where $w$ is an element of the Weyl group $W$ of $G$
and $t$ is an element of of a fixed maximal torus of the compact real form of $G$.
 
Let $\lambda \in P^+$, $V_\varepsilon(\lambda)$ be the irreducible
$U_\varepsilon({\mathfrak g})$-module of type ${\mathbf 1}$ with the
highest weight $\lambda$. Then $V_\varepsilon(\lambda)$ admits of a
positive definite hermitian form $(.,.)$, such that $(xv_1,v_2) = (v_1
x^*v_2)$ for all $v_1, v_2\in V_\varepsilon(\lambda)$, $x\in
U({\mathfrak g})$.

Let $\{ v^\nu_\mu\}$ be an orthogonal basis for weight space
$V_\varepsilon(\lambda)_\mu$, $\mu \in P^+$. Then $\bigcup \{v^\nu_\mu \}$
is an orthogonal basis for $V_\varepsilon(\lambda)$. Let
$C^\lambda_{\nu,s;\mu,r}(x):= (xv^r_\mu, v^s_\nu)$ be the associated
matrix element of $V_\varepsilon(\lambda)$. Then the matrix elements
$C^\lambda_{\nu,s;\mu,r}$ (where $\lambda$ runs through $P^+$, while
$(\mu,r)$ and $(\nu,s)$ runs independently through the index set of a 
basis of $V_\varepsilon(\lambda)$) form a basis of ${\mathcal
F}_\varepsilon(G)$ (see \cite{CP}).

Now every irreducible *-representations of ${\mathcal
F}_\varepsilon(\SL_2({\mathbf C}))$ is equivalent to a representation
belonging to one of the following two families, each of which is
parameterized by ${\mathbf S}^1 = \{ t\in {\mathbf C}\quad \vert \quad
\vert
t\vert = 1 \}$ -- (i) the family of one-dimensional representations
$\tau_t$, (ii) the family $\pi_t$ of representations in $\ell^2({\mathbf
N})$ (see \cite{CP}, page 437, Proposition 13.1.8).

Moreover, there exists a surjective homomorphism ${\mathcal
F}_\varepsilon(G) \to {\mathcal F}(\SL_2({\mathbf C}))$ induced by the
natural inclusion $\SL_2({\mathbf C}) \hookrightarrow G$ and by composing
the representation $\pi_{-1}$ of ${\mathcal F}(\SL_2({\mathbf C}))$ with
this homomorphism we obtain a representation of ${\mathcal
F}_\varepsilon(G)$ in $\ell^2({\mathbf N})$ denoted by $\pi_{s_i}$ where
$s_i$ appears in the reduced decomposition $w = s_{i_1} s_{i_2} \dots
s_{i_k}$. More precisely, $\pi_{s_i} : {\mathcal F}_\varepsilon(G) \to
{\mathcal L}(\ell^2({\mathbf N}))$ is of class CCR (see \cite{Di}), i.e.
its image is dense in the ideal of compact operators in ${\mathcal L} 
(\ell^2({\mathbf N}))$. 

The representation $\tau_t$ is one-dimensional and is of the form
$$\tau_t(C^\lambda_{\nu,s;\mu,r}) = \delta_{r,s}\delta_{\mu\nu}\exp(2\pi
\sqrt{-1}\mu(x)),$$ if 
$$t = \exp(2\pi\sqrt{-1}x)) \in {\mathbf T}, \mbox{ for } x \in {\mathfrak
t} 
= \Lie\;{\mathbf T},$$ see \S 13.1.D from \cite{CP}.

The unitary representations of compact quantum groups
are classified as follows.

\begin{itemize}
\item Every irreducible unitary representation of ${\mathcal
F}_\varepsilon(G)$ on a Hilbert space  is a completion of a unitarizable
highest weight representation. Moreover, two such representations are
equivalent if and only if they have the same weight, (\cite{CP}, Theorem
13.1.7).
\item The highest weight representations can be described as follows. Let 
$w = s_{i_1}.s_{i_2}. \dots s_{i_k}$ be a reduced decomposition of an
element $w$ of the Weyl group $W$. Then, (i) the Hilbert space tensor
product 
$$\rho_{w,t} = \pi_{s_{i_1}} \otimes \dots \otimes \pi_{s_{i_k}} \otimes
\tau_t$$ is an irreducible *-representation of ${\mathcal
F}_\varepsilon(G)$ which is associated to the Schubert cell ${\mathbf
S}_w$; (ii) up to equivalence, the representation $\rho_{w,t}$ does not
depend on the choice of the reduced decomposition of $w$; (iii) every
irreducible *-representation of ${\mathcal F}_\varepsilon(G)$ is
equivalent to some $\rho_{w,t}$, see (\cite{CP}, Theorem 13.1.9). 
\end{itemize}

Moreover, one can show that 
$$\bigcap_{(w,t)\in W \times {\mathbf T}}\ker \rho_{w,t} = \{ e \},$$ i.e.
the
representation 
$$\bigoplus_{w\in W} \int_{\mathbf T}^\oplus \rho_{w,t} dt$$ is a faithful
representation and 
$$\dim \rho_{w,t} = \left\{ \begin{array}{cl} 1, & \qquad \mbox{ if } w =
e,\\
				             \infty, & \qquad \mbox{ if } 
w \ne e
\end{array} \right. $$  (see Remarks [1] after Theorem 13.1.9
from\cite{CP}.)

We recall now the definition of C*-algebraic compact quantum group. 
Let $G$ be a complex algebraic group with compact real form.
\begin{definition}[\cite{CP}, Definition. 13.3.1] The {\sl C*-algebraic
compact
quantum group} $C^*_\varepsilon(G)$ is the C*-completion
 of the *-algebra ${\mathcal F}_\varepsilon(G)$ with respect to the
C*-norm $$\Vert f \Vert := \sup_\rho \Vert \rho(f) \Vert \quad (f\in
{\mathcal F}_\varepsilon(G)),$$ where $\rho$ runs through the
*-representations of ${\mathcal F}_\varepsilon(G)$ and the norm on the
right-hand side is the operator norm.  \end{definition}

We have the following result about the structure of compact quantum groups.
\begin{theorem} 
$$C^*_\varepsilon(G) \cong {\mathbf C}({\mathbf T}) \oplus \bigoplus_{e
\ne w\in
W}\int_{\mathbf T}^{\oplus} {\mathcal K}({\mathbf H}_{w,t}) dt,$$
where ${\mathbf C}({\mathbf T})$ is the algebra of continuous function on ${\mathbf
T}$ and ${\mathcal K}({\mathbf H})$ is the ideal of compact operators in a separable
Hilbert space ${\mathbf H}$ 
\end{theorem}
\begin{pf}
Let $w = s_{i_1}.s_{i_2}\dots s_{i_k}$ be a reduced decomposition of
the element $w\in W$ into a product of reflections. By construction, the
representations $\rho_{w,t} = \pi_{s_{i_1}} \otimes \dots \otimes
\pi_{s_{i_k}} \otimes \tau_t$ where $\pi_{s_i}$ is the composition of the
homomorphism of
${\mathcal F}_\varepsilon(G)$ onto ${\mathcal
F}_\varepsilon(\SL_2({\mathbf C}))$ and the representation $\pi_{-1}$ of
${\mathcal F}_\varepsilon(\SL_2({\mathbf C}))$  in the Hilbert space
$\ell^2({\mathbf N})$. Thus we have
$$\CD
\pi_{s_i} : C^*_\varepsilon(G) @>>> \> 
C^*_\varepsilon(\SL_2({\mathbf C}))
@>\pi_{-1}>> {\mathcal L}(\ell^2({\mathbf N}))
\endCD$$
The last one is $CCR$ (see Dixmier \cite{Di}) and we have
$$\pi_{s_i}(C^*_\varepsilon(G))\cong {\mathcal K}({\mathbf H}_{s_i}),
\quad \tau_{t}(C^*_\varepsilon(G))\cong {\mathbf C},$$
$$\rho_{w,t}(C^*_\varepsilon(G))\cong {\mathcal K}({\mathbf H}_{s_1})
\otimes \dots \otimes {\mathcal K}({\mathbf H}_{s_k})\otimes {\mathbf C} 
\cong {\mathcal
K}({\mathbf H}_{w,t}),$$
where ${\mathbf H}_{w,t} = {\mathbf H}_{s_1} \otimes \dots \otimes {\mathbf H}_{s_k}
\otimes {\mathbf C}.$ (Here the tensor product of Banach algebras always
means the tensor
with the minimal norm.) 
\end{pf}

\section{Stability of Entire Homology of Non-Commutative de Rham Currents}

We briefly recall in this subsection a definition of
operator
K-functors.  Following G. G. Kasparov, a
Fredholm representation of a C*-algebra $A$ is a triple
$(\pi_1,\pi_2,F)$,
consisting of *-representations $\pi_1, \pi_2 : A \to {\mathcal
L}({\mathcal H}_B)$ and a
Fredholm operator $F\in {\mathcal F}({\mathcal H}_B)$,
admitting an adjoint operator, on the
Hilbert
C*-module ${\mathcal H}_B= \ell^2_B$ over C*-algebra $B$, satisfying the
relations
$$\pi_1(a)F -
F\pi_2(a) \in {\mathcal K}_B,$$ where ${\mathcal K}_B$ is the ideal of compact
(adjoint-able) C*-module endomorphisms of ${\mathcal H}_B$.
The classes of homotopy invariance and unitary equivalence of Fredholm
modules form the
so called Kasparov
operator K-group $KK^*(A,B)$. Herewith put $A=C^*_\varepsilon(G)$, $B =
{\mathbf C}$,
where $G$ is a compact Lie group. In this case, we have
$K_*(C^*_\varepsilon(G)) =
KK^*(C^*_\varepsilon(G),{\mathbf C})$, where $K_*(A)$ is algebraic
$K$-group of $A$. 
From now on we shall use $K_*(A)$ to denote either topological or
 algebraic $K$-theory of the
Banach
algebra $A$. Note that
$$KK_0(A, {\mathbf C}) \cong K_0(A),$$
$$KK_1(A, {\mathbf C}) \cong K_1(A).$$

Let us recall now the definition and properties of entire homology of Banach
algebras. We refer the reader to \cite{DT1} for a more detailed exposition.
 
Let $A$ be an involutive Banach algebra. Recall that A. Connes defined entire cyclic
cohomology $\HE^*(A)$ and a pairing
$$K_*(A) \times \HE^*(A) \to {\mathbf C}.$$ Also M. Khalkhali \cite{Kh1}
proved
Morita and homotopy invariance of $\HE^*(A)$. We now define the entire
homology
$\HE_*(A)$ as follows: Given a collection $\{I_\alpha\}_{\alpha\in\Gamma}$
of ideals in $A$, equipped with a so called $\ad_A$-invariant trace
$$\tau_\alpha : I_\alpha
\to {\mathbf C},$$ satisfying the properties: \begin{enumerate}
\item $\tau_\alpha$ is a {\it continuous linear} functional, normalized as
$\Vert\tau_\alpha\Vert = 1$,
\item $\tau_\alpha$ is {\it positive} in the sense that $$\tau_\alpha(a^*a)
\geq 0, \forall
\alpha \in \Gamma,$$ where the map $a \mapsto a^*$ is the involution
defining the
involutive Banach algebra structure, i.e. an anti-hermitian endomorphism
such that $a^{**}
= a$ \item $\tau_\alpha$ is {\it strictly positive} in the sense that
$\tau_\alpha(a^*a) = 0$ iff
$a=0$, for every $\alpha\in \Gamma$. \item $\tau_\alpha$ is $\ad_A$-{\it
invariant} in the
sense that $$\tau_\alpha(xa) = \tau_\alpha(ax), \forall x\in A, 
a\in I_\alpha,$$
\end{enumerate}
then we have for every $\alpha \in \Gamma$ a scalar product $$\langle
a,b\rangle_\alpha
:=\tau_\alpha(a^*b)$$ and also an inverse system $\{I_\alpha,
\tau_\alpha\}_{\alpha \in
\Gamma}$. Let $\bar{I}_\alpha$ be the completion of $I_\alpha$ with respect to the
scalar product, defined
above and $\widetilde{\bar{I}_\alpha}$ denote $\bar{I}_\alpha$ with formally
adjoined unity
element. Define $C^n(\widetilde{\bar{I}_\alpha})$ as the set of
$n+1$-linear maps $\varphi :
(\widetilde{\bar{I}_\alpha})^{\otimes(n+1)} \to {\mathbf C}$. There exists a Hilbert
structure on  $(\widetilde{\bar{I}_\alpha})^{\otimes(n+1)}$ and we can identify
$C_n((\widetilde{\bar{I}_\alpha})) := \Hom(C^n(\bar{\tilde{I}}_\alpha), {\mathbf
C})$
with $C^n(\bar{\tilde{I}}_\alpha)$ via an anti-isomorphism.

For $I_\alpha \subseteq I_\beta$, we have a
well-defined map $$D^\beta_\alpha :
C^n(\widetilde{\bar{I}_\alpha}) \to C^n(\widetilde{\bar{I}_\beta}),$$ which
makes $\{
C^n(\widetilde{\bar{I}}_\alpha)\}$ into a direct system. Write $Q = \varinjlim
C^n(\widetilde{\bar{I}_\alpha})$ and observe that it admits a Hilbert
space 
structure, see \cite{DT1}-\cite{DT2}. Let $C_n(A) := \Hom(\varinjlim
C^n(\widetilde{\bar{I}_\alpha}), {\mathbf C}) = \Hom(\varinjlim
C_n(\widetilde{\bar{I}_\alpha}), {\mathbf C})$ which is anti-isomorphic to
$\varinjlim_\alpha C^n(\bar{\tilde{I}}_\alpha)$. So we have finally
$$C_n(A) = \varinjlim C_n(\bar{\tilde{I}}_\alpha).$$

Let $$b, b' : C^n(\widetilde{\bar{I}}_\alpha) \to
C^{n+1}(\widetilde{\bar{I}}_\alpha),$$
$$N : C^n(\widetilde{\bar{I}}_\alpha) \to
C^n(\widetilde{\bar{I}}_\alpha),$$ $$\lambda :
C^n(\widetilde{\bar{I}}_\alpha) \to C^n(\widetilde{\bar{I}}_\alpha),$$
$$S : C^{n+1}(\widetilde{\bar{I}}_\alpha)\to
C^n(\widetilde{\bar{I}}_\alpha)$$ be defined as in A. Connes \cite{Ca}.
We adopt the notations in \cite{Kh1} . Denote by $b^*, (b')^*, N^*,
\lambda^*, S^*$ the corresponding adjoint operators. Note also that for
each
$I_\alpha$ we have the same formulae for adjoint operators for homology as
Connes obtained for cohomology.

We now have a bi-complex
$$
\leqno{{\mathcal C}(A):}\qquad  
\CD
 @.    \vdots  @.    \vdots @.     \vdots @.      \\
@.          @V(-b')^* VV        @Vb^* VV       @V(-b')^* VV        @.   \\
\dots @<1-\lambda^* << \varinjlim_\alpha C_1(\bar{\tilde{I}}_\alpha) @<N^*<<
\varinjlim_\alpha C_1(\bar{\tilde{I}}_\alpha)
@<1-\lambda^*<< \varinjlim_\alpha C_1(\bar{\tilde{I}}_\alpha) @<N^*<< \dots\\
@.          @V(-b')^* VV        @Vb^* VV       @V(-b')^* VV        @.   \\
\dots @<1-\lambda^*<< \varinjlim_\alpha C_0(\bar{\tilde{I}}_\alpha) @<N^*<<
\varinjlim_\alpha 
C_0(\bar{\tilde{I}}_\alpha) @<1-\lambda^*<<\varinjlim_\alpha
C_0(\bar{\tilde{I}}_\alpha) @<N^*<< \dots\\
\endCD$$
with $d_v := b^*$ in the even columns and $d_v:= (-b')^*$ in the odd columns, $d_h
:= 1-\lambda^*$ from odd to even columns and $d_h := N^*$ from even to odd columns, 
where * means the corresponding adjoint operator. Now we have
$$\Tot({\mathcal C}(A))^{even} = \Tot({\mathcal C}(A))^{odd} :=
\oplus_{n\geq 0}
C_n(A),$$ which is periodic with period two. Hence, we have
$$\begin{array}{ccc}
&\partial & \\
\oplus_{n\geq 0} C_n(A)&\begin{array}{c} \longleftarrow\\ \longrightarrow
\end{array} &\oplus_{n\geq 0} C_n(A)\\
&\partial & \\
\end{array}$$
where $\partial = d_v + d_h$ is the total differential.

\begin{definition} Let $\HP_*(A)$ be the homology of the total complex
$(\Tot{\mathcal C}(A))$. It is called the {\sl periodic cyclic homology} 
of $A$.
\end{definition}

Note that this $\HP_*(A)$ is, in general different from the 
$\HP_*(A)$ of Cuntz-Quillen, because we used the direct limit of periodic
cyclic homology of ideals. But in special cases, when the whole
algebra $A$ is one of these ideals
 with $\ad_A$-invariant trace, (e.g. the
commutative algebras
of complex-valued functions on compact spaces), we return to the
Cuntz-Quillen $\HP_*$, which we shall use later. 

\begin{definition}
An even (or odd) chain $(f_n)_{n\geq 0}$ in ${\mathcal C}(A)$ is called
{\it entire} if the
radius of convergence of the power series $\sum_n \frac{n!}{[\frac{n}{2}]!}
\Vert
f_n\Vert z^n$, $z\in {\mathbf C}$ is infinite.
\end{definition}

Let $C^e(A)$ be the sub-complex of $C(A)$ consisting of entire chains. Then
we have a
periodic complex.

\begin{theorem} Let
$$\Tot(C^e(A))^{even} = \Tot(C_e(A))^{odd} := \bigoplus_{n\geq 0} C^e_n(A),$$
where $C^e_n(A)$ is the entire $n$-chain. Then we have a complex of
entire chains with the total differential $\partial := d_v + d_h$
$$\begin{array}{ccc} & \partial
& \\
\bigoplus_{n\geq 0} C^e_n(A) & \begin{array}{c}\longleftarrow \\
\longrightarrow\end{array} &
\bigoplus_{n\geq 0}
C^e_n(A)\\
& \partial & \end{array}$$
\end{theorem}
\begin{definition}
The homology of this complex is called also the {\it entire homology} and
denoted by $\HE_*(A)$. 
\end{definition}
Note that this entire homology is defined through
the inductive limits of ideals with ad-invariance trace. 

In \cite{DT1}-\cite{DT2}, the main properties of this theory, namely
\begin{itemize}
\item Homotopy invariance,
\item Morita invariance and
\item Excision,
\end{itemize}
were proved and hence that $\HE_*$ is a generalized homology theory.

\begin{lemma}
If the Banach algebra $A$ can be presented as a direct limit
$\varinjlim_\alpha I_\alpha$ of a system
of ideals $I_\alpha$ with $\ad$-invariant trace $\tau_\alpha$, then
$$K_*(A) = \varinjlim_\alpha K_*(I_\alpha),$$
$$\HE_*(A) = \varinjlim_\alpha \HE_*(I_\alpha).$$
\end{lemma}
\begin{pf}
Let us consider $K_*(A) = K_0(A) \oplus K_1(A)$. By definition 
$$ K_0(A) := \{ e^* = e = e^2 \in \Mat_\infty(A) \},$$ 
where 
$$\Mat_\infty(A) = \varinjlim_n Mat_n(A) = \varinjlim_n
\Mat_n(\varinjlim_\alpha
I_\alpha) = \varinjlim_n \varinjlim_\alpha \Mat_n(I_\alpha) = $$
$$= \varinjlim_\alpha  \varinjlim_n \Mat_n(I_\alpha) = \varinjlim_\alpha
\Mat_\infty( I_\alpha).$$ 
This means that $K_0(A) = \varinjlim_\alpha(I_\alpha)$. By the same
argument
we can prove for the $K_1$ groups: 
$$\begin{array}{rl} 
K_1(A) &= \GL_\infty(A) /(\GL_\infty(A),\GL_\infty(A))\\
       &\cong \GL_\infty(\varinjlim_\alpha I_\alpha)
/(\GL_\infty(\varinjlim_\alpha
I_\alpha),\GL_\infty(\varinjlim_\alpha I_\alpha)) \\
       &\cong \varinjlim_\alpha 
\GL(I_\alpha) / (\GL(I_\alpha),\GL(I_\alpha)) \\
       &\cong \varinjlim_\alpha K_*(I_\alpha).
\end{array}$$
We consider now the ${\mathcal C}(\bar{\tilde{I}}_\alpha)$-bi-complex
$$
\leqno{{\mathcal C}(\bar{\tilde{I}}_\alpha):}\qquad 
\CD
 @.    \vdots  @.    \vdots @.     \vdots @.      \\
@.          @V(-b')^* VV        @Vb^* VV       @V(-b')^* VV        @.   \\
\dots @<1-\lambda^* << C_1(\bar{\tilde{I}}_\alpha) @<N^* <<
C_1(\bar{\tilde{I}}_\alpha)
@<1-\lambda^* << C_1(\bar{\tilde{I}}_\alpha) @<N^* << \dots\\
@.          @V(-b')^* VV        @Vb^* VV       @V(-b')^*VV        @.   \\
\dots @<1-\lambda^* << C_0(\bar{\tilde{I}}_\alpha) @<N^* <<
C_0(\bar{\tilde{I}}_\alpha)
@<1-\lambda^* << C_0(\bar{\tilde{I}}_\alpha) @<N^* << \dots
\endCD$$
The connecting maps of the direct system 
$$\CD 
\bar{\tilde{I}}_\alpha @>>> \bar{\tilde{I}}_\beta @>>>
\bar{\tilde{I}}_\gamma, \mbox{ for } \alpha \leq
\beta \leq \gamma
\endCD$$
induces the same connecting maps for complexes
$$\CD 
{\mathcal C}(\bar{\tilde{I}}_\alpha) @>>> {\mathcal C}(\bar{\tilde{I}}_\beta) 
@>>>   {\mathcal C}(\bar{\tilde{I}}_\gamma), \mbox{ for } \alpha \leq
\beta \leq \gamma
\endCD$$
and finally we have for the corresponding total complexes:
$$\varinjlim_\alpha \Tot{\mathcal C}(\bar{\tilde{I}}_\alpha) = \bigoplus_\alpha
\Tot{\mathcal C}(\bar{\tilde{I}}_\alpha) / \sim .$$
It is easy to check the following claims:
$${\mathcal C}( \oplus_\alpha \bar{\tilde{I}}_\alpha) \cong
\bigoplus_\alpha {\mathcal
C}(\bar{\tilde{I}}_\alpha),$$ 
$$\Tot{\mathcal C}( \oplus_\alpha \bar{\tilde{I}}_\alpha) \cong
\bigoplus_\alpha \Tot{\mathcal
C}(\bar{\tilde{I}}_\alpha),$$
and finally
$$\Tot{\mathcal C}(\varinjlim_\alpha \bar{\tilde{I}}_\alpha) \cong \varinjlim
\Tot{\mathcal C}(\bar{\tilde{I}}_\alpha).$$
From the last claim, it is easy to see that a cycle of the complex 
$\Tot{\mathcal C}(\varinjlim_\alpha \bar{\tilde{I}}_\alpha)$ can be  
presented as a sequence of cycles as an element from $\varinjlim
\Tot{\mathcal C}(\bar{\tilde{I}}_\alpha)$.
The lemma is therefore proved.
\end{pf}

The following result from K-theory is well-known:
\begin{theorem} The entire homology of non-commutative de Rham currents
admits the following stability property
$$K_*({\mathcal K}({\mathbf H})) \cong K_*({\mathbf C}),$$
$$K_*(A \otimes {\mathcal K}({\mathbf H})) \cong K_*(A),$$
where ${\mathbf H}$ is a separable Hilbert space and $A$ is an arbitrary Banach
space.
\end{theorem}

The similar result is true for entire homology $\HE_*$ :
\begin{theorem} The entire homology of non-commutative de Rham currents
admits the following stability property
$$\HE_*({\mathcal K}({\mathbf H})) \cong \HE_*({\mathbf C}),$$
$$\HE_*(A \otimes {\mathcal K}({\mathbf H})) \cong \HE_*(A),$$
where ${\mathbf H}$ is a separable Hilbert space and $A$ is an arbitrary Banach
space.
\end{theorem}
\begin{pf}
It is easy to see that ${\mathcal K}({\mathbf H})$ is a direct
limit of matrix algebras $\Mat_n({\mathbf C})$ as
$${\mathcal K}({\mathbf H}) = \varinjlim_n \Mat_n({\mathbf C}),$$
and thus
$$A \otimes {\mathcal K}({\mathbf H}) = A \otimes \varinjlim_n
\Mat_n({\mathbf C}) = \varinjlim_n A \otimes \Mat_n({\mathbf C}).$$
The theorem is therefore deduced from the previous lemma.
\end{pf}

\section{Non-Commutative Chern Characters}

Let $A$ be a involutive Banach  algebra. 
We construct a non-commutative character $$ch_{C^*} : K_*(A) \to \HE_*(A)$$ and
later show
that when $A= C^*_\varepsilon(G)$, this Chern character reduces up to isomorphism to
classical Chern character on the normalizers of maximal compact tori. 

Let $A$ be an involutive Banach algebra with unity. 
\begin{theorem}
There exists a Chern character $$ch_{C^*} : K_*(A) \to \HE_*(A).$$ 
\end{theorem}
\begin{pf}
We first recall that there exists a pairing $$K_n(A) \times
C^n(\bar{\tilde{I}}_\alpha) \to
{\mathbf C}$$ due to
A. Connes, see \cite{Co}. Hence there exists a map
$ K_n(A) \stackrel {C_n }{ \longrightarrow }\Hom(C^n(A),{\mathbf C}).$ So,
by 1.1,
we have for each $\alpha$, a map $K_n(A)
\stackrel{C_n^\alpha}{\longrightarrow}
\Hom(C^n(\widetilde{\bar{I}_\alpha}),{\mathbf C})$ and hence a map $K_n(A)
\stackrel{C_n}{\longrightarrow}
\Hom(\varinjlim_\alpha C^n(\widetilde{\bar{I}_\alpha}),{\mathbf C}).$ We
now show that $C_n$
induces
the Chern map
$$ch: K_n(A) \to \HE_n(A)$$

Now let $e$ be an idempotent in $M_k(A)$ for some $k\in {\mathbf N}$. It
suffices to
show that for $n$ even, if $\varphi = \partial\psi$, where $\varphi\in
C^n(\widetilde{\bar{I}_\alpha})$ and $\psi \in
C^{n+1}(\widetilde{\bar{I}_\alpha})$,
then $$\langle e,\varphi\rangle = \sum_{n=1}^\infty
\frac{(-1)^n}{n!}\varphi(e,e,\dots,e)
= 0.$$
However, this follows from Connes' results in (\cite{Co},Lemma 7).

The proof of the case for $n$ odd would also follow from \cite{Co}. \end{pf}

Our next result computes the Chern characters  for $A=C^*_\varepsilon(G)$ by
reducing it to the classical case.

\begin{theorem}
Let ${\mathbf T}$ be a fixed maximal torus of $G$ with Weyl group $W:=
{\mathcal N}_{\mathbf T}/{\mathbf T}$.
Then the Chern character $$ch_{C^*}: K_*(C^*_\varepsilon(G)) \to
\HE_*(C^*_\varepsilon(G))$$ is an
isomorphism modulo torsion, i.e.
$$\CD
ch_{C^*}: K_* (C^*_\varepsilon(G))\otimes {\mathbf C} @>\cong >>
\HE_*(C^*_\varepsilon(G)),\endCD $$
which can be identified with the classical
Chern character
$$\CD ch: K_*({\mathbf C}({\mathcal N}_{\mathbf T})) @>>> \HE_*({\mathbf
C}({\mathcal N}_{\mathbf T}))\endCD$$
that
is also an isomorphism modulo torsion, i.e.
$$\CD
ch: K_* ({\mathcal N}_{\mathbf T})\otimes {\mathbf C} @>\cong >>
H^*_{DR}({\mathcal N}_{\mathbf T}).
\endCD$$ 
\end{theorem}

To prove this theorem we first observe:  

\begin{lemma}
$$K_*(C^*_\varepsilon(G)) \cong K_*({\mathbf C}({\mathcal N}_{\mathbf
T})) \cong K^*({\mathcal N}_{\mathbf T}).$$ 
\end{lemma}
\begin{pf}  We have
$$\begin{array}{cl}  
K_*(C^*_\varepsilon(G)) &= K_*({\mathbf C}({\mathbf T}) \oplus \bigoplus_{e
\ne w\in W}\int_{\mathbf T}^{\oplus} {\mathcal K}({\mathbf H}_{w,t}) dt)\\
  &= K_*({\mathbf C}({\mathbf T}) \oplus K_*(\bigoplus_{e
\ne w\in W}\int_{\mathbf T}^{\oplus} {\mathcal K}({\mathbf H}_{w,t}) dt)\\
  &= K_*({\mathbf C}({\mathbf T})) \oplus K_*(\bigoplus_{e
\ne w\in W}\int_{\mathbf T}^{\oplus} \varinjlim_n \Mat_n({\mathbf C})
dt)\\
  &= K_*({\mathbf C}({\mathbf T})) \oplus K_*(\varinjlim_n\bigoplus_{e
\ne w\in W}\int_{\mathbf T}^{\oplus}  \Mat_n({\mathbf C}) dt)\\
  &=K_*({\mathbf C}({\mathbf T})) \oplus \varinjlim_n K_*({\mathbf C}((W-\{e\})\times {\mathbf T}){\otimes}
\Mat_n({\mathbf C}))\\
  &\cong K_*({\mathbf C}({\mathbf T})) \oplus K_*({\mathbf C}((W-\{e\})\times 
{\mathbf T}) )\\
  &= K_*({\mathbf C}({\mathbf T}) \oplus {\mathbf C}((W-\{e\})\times 
{\mathbf T}) )\\
  &= K_*({\mathbf C}(W \times {\mathbf T}) )\\ 
  &= K_*({\mathbf C}({\mathcal N}_{\mathbf T}/{\mathbf T} \times {\mathbf T}))\\
  &= K_*({\mathbf C}({\mathcal N}_{\mathbf T})),
\end{array}$$ 
by using Morita invariance property. The last one is isomorphic with the
topological K-theory $K^*({\mathcal N}_{\mathbf T})$. 
\end{pf}

Next we have 
\begin{lemma}
$$\HE_*(C^*_\varepsilon(G)) \cong \HE_*({\mathbf C}({\mathcal N}_{\mathbf
T})) \cong H^*_{DR}({\mathcal N}_{\mathbf T}).$$ 
\end{lemma}
\begin{pf} The proof is analogous to that of 3.3. So we have
$$\begin{array}{cl}  
\HE_*(C^*_\varepsilon(G)) &= \HE_*({\mathbf C}({\mathbf T}) \oplus \bigoplus_{e
\ne w\in W}\int_{\mathbf T}^{\oplus} {\mathcal K}({\mathbf H}_{w,t}) dt)\\
  &= \HE_*({\mathbf C}({\mathbf T}) \oplus \HE_*(\bigoplus_{e
\ne w\in W}\int_{\mathbf T}^{\oplus} {\mathcal K}({\mathbf H}_{w,t}) dt)\\
  &= \HE_*({\mathbf C}({\mathbf T})) \oplus \HE_*(\bigoplus_{e
\ne w\in W}\int_{\mathbf T}^{\oplus} \varinjlim_n \Mat_n({\mathbf C}) dt)\\
  &= \HE_*({\mathbf C}({\mathbf T})) \oplus \HE_*(\varinjlim_n\bigoplus_{e
\ne w\in W}\int_{\mathbf T}^{\oplus}  \Mat_n({\mathbf C}) dt)\\
  &=\HE_*({\mathbf C}({\mathbf T})) \oplus \varinjlim_n \HE_*({\mathbf C}((W-\{e\})\times {\mathbf T}){\otimes}
\Mat_n({\mathbf C}))\\
  &\cong \HE_*({\mathbf C}({\mathbf T})) \oplus \HE_*({\mathbf C}((W-\{e\})\times 
{\mathbf T}) )\\
  &= \HE_*({\mathbf C}({\mathbf T}) \oplus {\mathbf C}((W-\{e\})\times 
{\mathbf T}) )\\
  &= \HE_*({\mathbf C}(W \times {\mathbf T}) )\\ 
  &= \HE_*({\mathbf C}({\mathcal N}_{\mathbf T}/{\mathbf T} \times {\mathbf T}))\\
  &= \HE_*({\mathbf C}({\mathcal N}_{\mathbf T})),
\end{array}$$ 
by using Morita invariance property of $\HE_*$. Now, since ${\mathbf
C}({\mathcal N}_{\mathbf T})$ is a commutative ${\mathbf C}$-algebra, then
by Cuntz-Quillen's result \cite{CQ}, we have a canonical isomorphism
$\HP_*({\mathbf C}({\mathcal N}_{\mathbf T})) \cong H^*_{DR}({\mathcal
N}_{\mathbf T})$, since ${\mathbf C}({\mathcal N}_{\mathbf T})$ is smooth.
Moreover by a result of Khalkhali \cite{Kh1}-\cite{Kh2} $\HP_*({\mathbf
C}({\mathcal N}_{\mathbf T})) \cong 
\HE_*({\mathbf C}({\mathcal N}_{\mathbf T}))$; hence the result.
\end{pf}

{\sl Proof of Theorem 3.2}

Now consider the commutative diagram
$$\leqno{(I)} \qquad\qquad\qquad
\CD
K_* (C^*_\varepsilon(G)) @>ch_{C^*}>> \HE_*(C^*_\varepsilon(G))\\
@V\alpha VV    @VV\gamma V\\
K_* ({\mathbf C}({\mathcal N}_{\mathbf T})) @>ch_{CQ}>>
\HE_*({\mathbf C}({\mathcal N}_{\mathbf T}))\\
@V\beta VV    @VV\delta V\\
K^* ({\mathcal N}_{\mathbf T}) @>ch>> H^*_{DR}({\mathcal
N}_{\mathbf T})
\endCD$$ 

Now, $K_* (C^*_\varepsilon(G)) \cong K_* ({\mathbf C}({\mathcal
N}_{\mathbf T}))$ by Lemma 3.3. Also $\HE_*(C^*_\varepsilon(G)) \cong
\HE_*({\mathbf C}({\mathcal N}_{\mathbf T}))$ by Lemma 3.4. By a result of
Khalkhali \cite{Kh1}-\cite{Kh2}, we have $\HE_*({\mathbf C}({\mathcal
N}_{\mathbf T})) \cong \HP_*({\mathbf C}({\mathcal N}_{\mathbf T}))$.
Moreover $\HP_*({\mathbf C}({\mathcal N}_{\mathbf T})) \cong
H^*_{DR}({\mathcal N}_{\mathbf T})$, since the algebra ${\mathbf
C}({\mathcal N}_{\mathbf T})$ is smooth in the sense of Cuntz-Quillen (see
\cite{CQ}). Hence, $\HE_*({\mathbf C}({\mathcal N}_{\mathbf T})) \cong
H^*_{DR}({\mathcal N}_{\mathbf T})$. Furthermore, $K_*({\mathbf
C}({\mathcal N}_{\mathbf T})) \cong
K^*({\mathcal N}_{\mathbf T})$ by the well-known Swan' theorem.  
So all the vertical maps are isomorphisms. It is well known that the
bottom map $ch$ is
also an isomorphism modulo torsion, i.e.
$$\CD ch: K^* ({\mathcal N}_{\mathbf T})\otimes {\mathbf C} @>\cong >>
H^*_{DR}({\mathcal
N}_{\mathbf T}). \endCD$$
Hence $ch_{CQ}$ and
$ch_{C^*}$ are isomorphisms modulo torsion.
The theorem is therefore proved.

\section{Algebraic Version}

Let $C^*_\varepsilon(G)$ be a compact quantum group, $\HP_*(C^*_\varepsilon(G))$ the
periodic cyclic homology introduced in \S2. Since $C^*_\varepsilon(G) = {\mathbf
C}({\mathbf T}) \oplus \bigoplus_{e \ne w\in W}\int_{\mathbf T}^{\oplus} {\mathcal
K}({\mathbf H}_{w,t}) dt$

$\HP_*(C^*_\varepsilon(G))$ coincides with the $\HP_*(C^*_\varepsilon(G))$ defined
by J.  Cuntz-D.  Quillen \cite{CQ}.

J. Cuntz and D. Quillen \cite{CQ} defined the so called $X$-complexes of
${\mathbf C}$-algebras and then used some ideas of Fedosov product to
define algebraic Chern
characters. We now briefly recall their definitions. For a
(non-commutative) associate ${\mathbf C}$-algebra $A$, consider the space
of even
non-commutative differential forms $\Omega^+(A) \cong RA$, equipped with
the Fedosov
product
$$\omega_1 \circ \omega_2 := \omega_1\omega_2 - (-1)^{\vert \omega_1\vert}
d\omega_1
d\omega_2,$$ see \cite{CQ}. Consider also the ideal $IA := \oplus_{k\geq 1}
\Omega^{2k}(A)$. It is easy to see that $RA/IA \cong A$ and that $RA$
admits the
universal property that any based linear map $\rho : A \to M$ can be
uniquely extended to a
derivation $D : RA \to M$. The derivations $D : RA \to M$ bijectively
correspond to lifting
homomorphisms from $RA$ to the semi-direct product $RA \oplus M$, which also
bijectively correspond to linear map $\bar\rho : \bar{A}= A /{\mathbf C}
\to M$ given by $$
a\in \bar{A} \mapsto D(\rho a).$$ From the universal property of
$\Omega^1(RA)$, we
obtain a bi-module isomorphism $$RA \otimes \bar{A} \otimes RA \cong
\Omega^1(RA).$$ As in \cite{CQ}, let $\Omega^-A = \oplus_{k\geq 0}
\Omega^{2k+1}A$. Then we have $$\Omega^{-}A \cong RA \otimes \bar{A} \cong
\Omega^1(RA)_\# := \Omega^1(RA)/[(\Omega^1(RA),RA)].$$
\par
J. Cuntz and D. Quillen proved
\begin{theorem}(\cite{CQ}, Theorem1):
There exists an isomorphism of
${\mathbf Z}/(2)$-graded complexes
$$\Phi : \Omega A = \Omega^+A \oplus \Omega^{-}A \cong RA \oplus
\Omega^1(RA)_\#,$$ such that
$$\Phi : \Omega^+A \cong RA,$$ is defined by $$\Phi(a_0da_1\dots da_{2n} =
\rho(a_1)\omega(a_1,a_2) \dots \omega(a_{2n-1},a_{2n}),$$
and $$ \Phi : \Omega^{-}A \cong \Omega^1(RA)_\#,$$ $$\Phi(a_0da_1\dots
da_{2n+1})
= \rho(a_1)\omega(a_1,a_2)\dots \omega(a_{2n-1},a_{2n})\delta(a_{2n+1}).$$
With
respect to this identification, the product in $RA$ is just the Fedosov
product on even
differential forms and the differentials on the $X$-complex
$$X(RA) : \qquad RA\cong \Omega^+A \to \Omega^1(RA)_\# \cong \Omega^{-}A
\to RA
$$ become the operators
$$\beta = b - (1+\kappa)d : \Omega^{-}A \to \Omega^+A,$$ $$\delta =
-N_{\kappa^2} b
+ B : \Omega^+A \to \Omega^{-}A,$$ where $N_{\kappa^2} =
\sum_{j=0}^{n-1} \kappa^{2j}$, $\kappa(da_1\dots da_n) := da_n\dots da_1$.
\end{theorem}
\par
Let us denote by $IA \triangleleft RA$ the ideal of even non-commutative
differential forms
of
order $\geq 2$. By the universal property of $\Omega^1$ $$\Omega^1(RA/IA) =
\Omega^1RA/((IA)\Omega^1RA + \Omega^1RA.(IA) + dIA).$$ Since $\Omega^1RA =
(RA)dRA = dRA.(RA)$, then $\Omega^1RA(IA) \cong IA\Omega^1RA \mod
[RA,\Omega^1R].$
$$\Omega^1(RA/IA)_\# = \Omega^1RA /([RA,\Omega^1RA]+IA.dRA + dIA).$$ For
$IA$-adic tower $RA/(IA)^{n+1}$, we have the complex $${\mathcal
X}(RA/(IA)^{n+1}) : \qquad RA/IA^{n+1} \leftarrow
\Omega^1RA/([RA,\Omega^1RA]+(IA)^{n+1}dRA + d(IA)^{n+1}).$$
Define
$${\mathcal X}^{2n+1}(RA,IA) : \quad RA/(IA)^{n+1} \to
\Omega^1RA/([RA,\Omega^1RA]+(IA)^{n+1}dRA + d(IA)^{n+1}) $$ $$\to
RA/(IA)^{n+1},$$
$${\mathcal X}^{2n}(RA,IA): \quad RA/((IA)^{n+1} +[RA,IA^n]) \to
\Omega^1RA/([RA,\Omega^1RA]+d(IA)^ndRA)$$ $$\to RA/((IA)^{n+1}
+[RA,IA^n]).$$
One has
$$b((IA)^ndIA) = [(IA)^n,IA] \subset (IA)^{n+1},$$ $$d(IA)^{n+1} \subset
\sum_{j=0}^n (IA)^jd(IA)(IA)^{n-j} \subset (IA)^n dIA
+ [RA,\Omega^1RA].$$
and hence
$${\mathcal X}^1(RA,IA = X(RA,IA),$$
$${\mathcal X}^0(RA,IA) = (RA/IA)_\#.$$
There is a sequence of maps between complexes $$\dots \to X(RA/IA) \to
{\mathcal
X}^{2n+1}(RA,IA)\to {\mathcal X}^{2n}(RA,IA) \to X(RA/IA) \to \dots $$ We
have the
inverse limits
$$\hat{X}(RA,IA) := \varprojlim X(RA/(IA)^{n+1}) = \varprojlim {\mathcal
X}^n(RA,IA).$$
Remark that
$${\mathcal X}^q = \Omega A/F^q\Omega A,$$ $$\hat{X}(RA/IA)= \hat{\Omega}A.$$

We quote the second main result of J. Cuntz and D. Quillen (\cite{CQ},
Thm2), namely:

$$H_i\hat{\mathcal X}(RA,IA) = \HP_i(A).$$

We now apply  this machinery to our case. First we have the following algebraic
analogy of the Lemma for $\HP_*$.

\begin{lemma}
$$\HP_*(C^*_\varepsilon(G)) \cong \HP_*({\mathbf C}({\mathcal N}_{\mathbf
T})) \cong H^*_{DR}({\mathcal N}_{\mathbf T}).$$ 
\end{lemma}
\begin{pf} By analogy with the proof of the previous proposition, we also have
$$\begin{array}{cl}  
\HP_*(C^*_\varepsilon(G)) &= \HP_*({\mathbf C}({\mathbf T}) \oplus \bigoplus_{e
\ne w\in W}\int_{\mathbf T}^{\oplus} {\mathcal K}({\mathbf H}_{w,t}) dt)\\
  &= \HP_*({\mathbf C}({\mathbf T}) \oplus \HP_*(\bigoplus_{e
\ne w\in W}\int_{\mathbf T}^{\oplus} {\mathcal K}({\mathbf H}_{w,t}) dt)\\
  &= \HP_*({\mathbf C}({\mathbf T})) \oplus \HP_*(\bigoplus_{e
\ne w\in W}\int_{\mathbf T}^{\oplus} \varinjlim_n \Mat_n({\mathbf C}) dt)\\
  &= \HP_*({\mathbf C}({\mathbf T})) \oplus \HP_*(\varinjlim_n\bigoplus_{e
\ne w\in W}\int_{\mathbf T}^{\oplus}  \Mat_n({\mathbf C}) dt)\\
  &=\HP_*({\mathbf C}({\mathbf T})) \oplus \varinjlim_n \HP_*({\mathbf
C}((W-\{e\})\times {\mathbf T}){\otimes}
\Mat_n({\mathbf C}))\\
  &\cong \HP_*({\mathbf C}({\mathbf T})) \oplus \HP_*({\mathbf C}((W-\{e\})\times 
{\mathbf T}) )\\
  &= \HP_*({\mathbf C}({\mathbf T}) \oplus {\mathbf C}((W-\{e\})\times 
{\mathbf T}) )\\
  &= \HP_*({\mathbf C}(W \times {\mathbf T}) )\\ 
  &= \HP_*({\mathbf C}({\mathcal N}_{\mathbf T}/{\mathbf T} \times
{\mathbf T}))\\
  &= \HP_*({\mathbf C}({\mathcal N}_{\mathbf T})),
\end{array}$$ 
by using Morita invariance property of $\HP_*$. The last one is following
Cuntz-Quillen for the commutative ${\mathbf
C}$-algebra
$A$, we have a canonical isomorphism from periodic cyclic homology to the
${\mathbf
Z}/(2)$-graded de Rham homology, when $A$ is
smooth.
\end{pf}

Now, for each idempotent $e\in M_n(A)$ there is an unique element $x\in
M_n(\widehat{RA})$.
Then the element $$\tilde{e} := x + (x-\frac{1}{2})\sum_{n\geq 1} \frac{2^n(2n-
1)!!}{n!}(x-x^{2n})^{2n}\in M_n(\widehat{RA})$$ is a lifting of $e$ to an
idempotent
matrix in $M_n(\widehat{RA})$. Then the map $[e] \mapsto tr(\tilde{e})$
defines the map
$K_0(A) \to H_0(X(\widehat{RA})) = \HP_0(A)$. To an element $g\in\GL_n(A)$
one
associates an element $p\in \GL(\widehat{RA})$ and to the element $g^{-1}$
an element
$q\in \GL_n(\widehat{RA})$ then put
$$x = 1- qp, \mbox{ and } y = 1-pq.$$
And finally, to each class
$[g]\in \GL_n(A)$ one associates $$tr(g^{-1}dg) = tr(1-x)^{-1}d(1-x) =
d(tr(log(1-x))) =
-tr\sum_{n=0}^\infty x^ndx\in \Omega^1(A)_\#.$$ Then $[g] \to tr(g^{-1}dg)$
defines
the map $K_1(A) \to HH_1(A) = H_1(X(\widehat{RA})) = \HP_1(A)$.

\begin{definition} Let $\HP^*(I_n)$ be the periodic cyclic cohomology defined by
Cuntz-Quillen. Then the pairing
$$K_*^{alg}(C^*(G)) \times \bigcup_n \HP^*(I_n) \to {\mathbf C},$$ 
where $I_n = {\mathbf C}({\mathbf T}) \oplus \bigoplus \int^\oplus_{\mathbf T}
\Mat_n({\mathbf C}) \cong {\mathbf C}({\mathbf T}) \oplus{\mathbf
C}((W-\{e\})
\times {\mathbf T}) \otimes \Mat_n({\mathbf C}),$
 defines an
algebraic non-commutative Chern character $$ch_{alg} : K_*^{alg}(C^*_\varepsilon(G))
\to
\HP_*(C^*_\varepsilon(G)),$$ which gives us a variant of non-commutative Chern
characters with
values in $\HP_*$-groups. \end{definition}

We close this section with an algebraic analog of theorem 3.2.

\begin{theorem}
Let $G$ be a compact Lie group and ${\mathbf T}$ a fixed maximal
compact torus
of $G$.
Then in the notations of 4.3, the Chern character $$ch_{alg} :
K_*(C^*_\varepsilon(G)) \to
\HP_*(C^*_\varepsilon(G))$$ is an isomorphism modulo torsion,i.e.
$$ch_{alg} :
K_*(C^*_\varepsilon(G)) \otimes {\mathbf C}\to
\HP_*(C^*_\varepsilon(G))$$
 which
can be
identified with the
classical Chern,
$$\CD ch: K^* ({\mathcal N}_{\mathbf T}) @>>> H^*_{DR}({\mathcal
N}_{\mathbf T}),
\endCD$$
 which  is equivalent to the
Cuntz-Quillen 
character $$ch_{CQ}: K_*({\mathbf C}({\mathcal N}_{\mathbf T})) \to \HP_*({\mathbf
C}({\mathcal N}_{\mathbf T}))$$ 
which are also isomorphisms modulo torsion, i.e.
$$\CD ch: K^* ({\mathcal N}_{\mathbf T})\otimes {\mathbf C} @>\cong >>
H^*_{DR}({\mathcal
N}_{\mathbf T}),
\endCD$$
$$\CD ch_{CQ}: K_*({\mathbf C}({\mathcal N}_{\mathbf T}))\otimes {\mathbf
C} 
@>\cong >> \HP_*({\mathbf
C}({\mathcal N}_{\mathbf T})).\endCD$$ 
\end{theorem}
\begin{pf}

First note that 
$$K_*(C^*_\varepsilon(G)) \cong K_*({\mathbf C}({\mathcal N}_{\mathbf
T})) \cong K^*({\mathcal N}_{\mathbf T}).$$ 

Next we have
$$\HP_*(C^*_\varepsilon(G)) \cong \HP_*({\mathbf C}({\mathcal N}_{\mathbf
T})) \cong H^*_{DR}({\mathcal N}_{\mathbf T}).$$ 

Furthermore, by a result of Cuntz-Quillen for the commutative ${\mathbf
C}$-algebra
$A$, we have a canonical isomorphism from periodic cyclic homology to the
${\mathbf
Z}/(2)$-graded de Rham cohomology, when $A$ is
smooth.
Hence
$$K^*({\mathbf C}({\mathcal N}_{\mathbf T})) \cong H^*({\mathbf C}({\mathcal 
N}_{\mathbf T})).$$
Now we
have a commutative diagram
$$\CD
K_*( C^*_\varepsilon(G)) @>ch_{alg}>> \HP_*(C^*_\varepsilon(G))\\
@V\cong VV    @VV\cong V\\
K_* ({\mathbf C}({\mathcal N}_{\mathbf T})) @>ch_{CQ}>>
\HP_*({\mathbf C}({\mathcal N}_{\mathbf T}))\\
@V\cong VV    @VV\cong V\\
K^* ({\mathcal N}_{\mathbf T}) @>ch>> H^*_{DR}({\mathcal
N}_{\mathbf T})
\endCD$$ 
where the vertical homomorphisms all are isomorphisms. Hence, $ch_{alg}$ is an
isomorphism modulo torsion.
\end{pf}

\section*{Acknowledgments} This work was completed during the stay of the
first author as a visiting mathematician at the International Centre for
Theoretical Physics, Trieste, Italy.  He would like to thank ICTP for the
hospitality, without which this work would not have been possible and the
IMU Commission of Development and Exchange for a generous grant.

This work is supported in part by the International Centre for Theoretical
Physics, Trieste, Italy, the IMU Commission on Development and Exchange 
and the Vietnam National Foundation for Research in Fundamental Sciences.

\end{document}